\documentclass[12pt,twoside]{article}
\usepackage{amsmath,amssymb,amsthm}

\theoremstyle{plain} 
\newtheorem{thm}{Theorem}[section]

\newtheorem{lem}{Lemma}[section]
\theoremstyle{remark}
\newtheorem{rem}{Remark}[section]

\newtheorem{defi}{Definition}[section]
\newtheorem*{pf}{Proof}

\def\ee{\mathbb {E}}
\def\dd{\mathbb {D}}
\def\ii{\mathbb {I}}

\def\pp{\mathbb {P}}

\def\cov{{\rm cov}}
\def\var{{\rm var}}

\def\({\left(}
\def\){\right)}
\def\[{\left[}
\def\]{\right]}
\def\rom#1{{\rm{#1}}}
\def\E{{\mathbb E}}
\let\ee=\E

\begin{document}

\baselineskip 5.2truemm

\vskip0.3truecm
\centerline{\Large{\bf On a General Approach to the}}
\smallskip
\centerline{ \renewcommand{\thefootnote}{1}
\Large{\bf Strong Laws of Large Numbers}\footnotemark
\footnotetext{Partially supported by the Hungarian Foundation of Scientific Researches
  under Grant No. OTKA T047067/2004  and Grant No. OTKA T048544/2005.} }
\vskip0.3truecm

\centerline{
  \renewcommand{\thefootnote}{2}
  {\bf Istv\'an Fazekas}\footnotemark
  \footnotetext{Faculty of Informatics, University of Debrecen,  P.O. Box 12, 4010 Debrecen, Hungary,  e-mail: fazekasi@inf.unideb.hu,
  tel: 36-52-512900/75211}}

\medskip
\begin{abstract}
\noindent
A general method to obtain strong laws of large numbers is studied.
The method is based on abstract H\'ajek-R\'enyi type maximal inequalities.
The rate of convergence in the law of large numbers is also considered.
Some applications for weakly dependent sequences are given.
\end{abstract}



\section{Introduction}
\setcounter{equation}{0}
%
H\'ajek and R\'enyi (1955) proved the following inequality.

Let $X_1, \dots , X_n$ be independent random variables with zero mean values and finite variances $\ee X_k^2 = \sigma_k^2$, $k=1, \dots , n$.
Denote by $S_k = X_1 + \dots + X_k$, $k=1, \dots , n$, the partial sums.
Let $\beta_1, \dots,\beta_n$ be a non-decreasing sequence of positive numbers.  
Then for any $\varepsilon> 0$ and for any $m$ with $1\leq m \leq n$ we have
\begin{equation}
\label{HaRe}
\pp \left( \max_{m \leq l\leq n} \left|\frac{S_l}{\beta_l}\right| \ge \varepsilon \right)
\leq
 \frac{1}{\varepsilon^2} \left[ \frac{1}{\beta_m^2}\sum_{l=1}^m \sigma_l^2 +  \sum_{l=m+1}^n \frac{\sigma_l^2}{\beta_l^2} \right].
\end{equation}
In H\'ajek-R\'enyi (1955) this inequality was used to obtain strong laws of large numbers (SLLN).

H\'ajek-R\'enyi type inequalities were proved for arbitrary random variables by Kounias and Weng (1969) and by Szynal (1973),
for independent random variables by Bickel (1970), for martingales by Chow (1960).
Later on H\'ajek-R\'enyi type maximal inequalities were obtained for partial sums of certain dependent random variables.
After a H\'ajek-R\'enyi type  inequality is obtained the proof of the strong law of large numbers becomes a standard task.

In Fazekas and Klesov (2000) it was shown that a H\'ajek-R\'enyi type maximal inequality for moments is always 
a consequence of an appropriate Kolmogorov type maximal inequality.
Moreover, the  H\'ajek-R\'enyi type maximal inequality automatically implies the strong law of large numbers. 
The most important is that no restriction is assumed on the dependence structure of the random variables.

In T\'om\'acs-L\'{i}bor (2006) an abstract  H\'ajek-R\'enyi type inequality is presented for the probability and 
it is shown that a Kolmogorov type inequality always implies an appropriate SLLN.
The framework is similar to the one in Fazekas-Klesov (2000).

The methods used in Fazekas-Klesov (2000) and T\'om\'acs-L\'{i}bor (2006) are known for the specialists.
However, it is worth to present the general forms of the SLLN.
Those explicit forms enable us to give unified proofs for the existing SLLN's, to improve them, and to find new ones.
Moreover, using the general approach, some results on the rate of convergence in the SLLN can be obtained (Hu-Hu (2006)).

A survey of the known versions of the Kolmogorov and H\'ajek-R\'enyi inequalities shows that most of them can be inserted into a few standard classes.
In the present paper we list some abstract versions of the  Kolmogorov and H\'ajek-R\'enyi inequalities.
Then we find their general relationships.
We shall see that a H\'ajek-R\'enyi inequality is always a consequence of the appropriate Kolmogorov inequality.
Then we shall show that they automatically imply strong laws of large numbers.

We shall show that several  H\'ajek-R\'enyi type inequalities and SLLN's can be inserted into the framework of our general theory.
We shall mention results e.g. on martingales, mixingales, $\varrho$-mixing sequences, sequences with superadditive momentum structure,
associated sequences, (asymptotically almost) negatively associated sequences, and demimartingales.

The original version of the present paper was presented at Seminar on Stability Problems for Stochastic Models, Sovata-Bay (Romania), 2006.
The manuscript was presented on the author's home page (see Fazekas (2006)).
However, the results became known for certain group of specialist (e.g. cited by Wang-Hu-Shen-Ling (2008) and Wang-Hu-Shen-Yang (2010)).
Now we supplement the original version with a new section (Section \ref{uj}) and update the list of references.

We shall use the following notation. $X_1, X_2, \dots$, will denote a sequence of random variables defined on a fixed probability space. 
The partial sums of the random variables will be $S_n= \sum_{i=1}^n X_i$ for $n\geq 1$ and $S_0=0$.
A sequence $\{ b_n \}$ will be called non-decreasing if $b_i \leq b_{i+1}$ for $i\geq 1$.
%
\section{H\'ajek-R\'enyi type maximal inequalities for moments}
\setcounter{equation}{0}
%
We start with the abstract Kolmogorov type and H\'ajek-R\'enyi type maximal inequalities for moments.
\begin{defi}   \label{defKolm1m}
We say that the  random variables $X_1, \dots, X_n$  satisfy the first Kolmogorov type maximal inequality for moments,
if for each $m$ with $1\leq m \leq n$ 
\begin{equation}
\label{Kolm1m}
\ee \Big[ \max_{1\leq l\leq m} |S_l| \Big]^r
\leq
K \sum\nolimits_{l=1}^m \alpha_l
\end{equation}
where $\alpha_1, \dots,\alpha_n$ are non-negative numbers, $r>0$, and $K>0$.
\end{defi}

\begin{defi}   \label{defHaRe1m}
We say that the  random variables $X_1, \dots, X_n$  satisfy the first H\'ajek-R\'enyi type maximal inequality for moments,
if
\begin{equation}
\label{HaRe1m}
\ee \left[ \max_{1\leq l\leq n} \left|\frac{S_l}{\beta_l}\right| \right]^r
\leq
C \sum_{l=1}^n \frac{\alpha_l}{\beta_l^r} 
\end{equation}
where $\beta_1, \dots,\beta_n$ is a non-decreasing sequence of positive numbers,  
$\alpha_1, \dots,\alpha_n$ are non-negative numbers, 
$r>0$, and $C>0$.
\end{defi}

\begin{thm}
\label{Kolm1m-HaRe1m}
(Fazekas-Klesov (2000), Theorem 1.1.)
Let the  random variables $X_1, \dots, X_n$ be fixed.
If the first Kolmogorov type maximal inequality \eqref{Kolm1m} for moments is satisfied, 
then the first H\'ajek-R\'enyi type maximal inequality \eqref{HaRe1m} for moments is satisfied with $C= 4K$.
\end{thm}
Here and in what follows we mean that the H\'ajek-R\'enyi type maximal inequality is valid with the same parameters $n$, $\alpha_1, \dots,\alpha_n$, $r>0$ as the 
appropriate Kolmogorov's inequality and for arbitrary non-decreasing sequences $\beta_1, \dots,\beta_n$ of positive numbers.
Moreover, $C$ may depend on $K$ and $r$ only.
 
The above result allows us to obtain an abstract form of the strong law of large numbers.
\begin{thm}
\label{Kolm1m-SLLN} (Theorem 2.1 in Fazekas-Klesov (2000).)
Let $X_1, X_2, \dots$ be a sequence of random variables.
Let the non-negative numbers $\alpha_1, \alpha_2, \dots$, $r>0$, and $K>0$ be fixed.
Assume that for each $m\geq 1$ the first Kolmogorov type maximal inequality \eqref{Kolm1m} for moments is satisfied.
Let $b_1, b_2, \dots$ be a non-decreasing unbounded sequence of positive numbers. 
If
\begin{equation}
\label{sum}
\sum\nolimits_{l=1}^{\infty} \frac{\alpha_l}{b_l^r} < \infty ,
\end{equation}
then
\begin{equation}
\label{lim}
\lim_{n\to\infty} \frac{S_n}{b_n} = 0 \ \ \text {a.s.}
\end{equation}
\end{thm}

\begin{pf}
The proof in Fazekas-Klesov (2000) is based on a theorem of Dini (see Lemma \ref{LemmaDini} below).
Actually it suffices to apply Lemma~\ref{Lemma22} below which is a simple consequence of Dini's theorem.
Let $\{\beta_n\}$ be a sequence satisfying the properties given in Lemma~\ref{Lemma22}.
Then apply Theorem~\ref{Kolm1m-HaRe1m}.
\qed\end{pf}

\begin{lem}
\label{Lemma22}
Let $\{b_k\}$ be a non-decreasing unbounded sequence of positive numbers.
Let $\{\alpha_k\}$ be a sequence of non-negative numbers, with
$
\sum\nolimits_{k=1}^{\infty} \frac{\alpha_k}{b_k^r} < \infty ,
$
where $r>0$. 
Then there exists a non-decreasing unbounded sequence $\{\beta_k\}$ of positive numbers such that
$\sum_{k=1}^{\infty} \frac{\alpha_k}{\beta_k^r} < \infty$ and
$\lim_{k\to\infty} \frac{\beta_k}{b_k} = 0$.
\end{lem}
In Fazekas-Klesov (1998) a direct proof of Lemma~\ref{Lemma22} is given.

The following remark in some cases makes easier the application of Theorem~\ref{Kolm1m-SLLN}.
\begin{rem}
\label{Remark21} (Corollary 2.1 and Remark 2.1 in Fazekas-Klesov (2000).)
Let $b_1, b_2, \dots$ be a non-decreasing unbounded sequence of positive numbers. 
Let $\alpha_1, \alpha_2, \dots$ be non-negative numbers, $\Lambda_k=\alpha_1 + \cdots+\alpha_k$ for $k\geq 1$.
Let $r$ be a fixed positive number. 
If
\begin{equation}
\label{eq24}
\sum\nolimits_{l=1}^{\infty} \Lambda_l\left({1}/{b_l^r} - {1}/{b_{l+1}^r}\right)
<\infty
\end{equation}
and
\begin{equation}
\label{eq25}
{\Lambda_n}/{b_n^r} \qquad \text{is bounded,}
\end{equation}
then (\ref{sum}) is fulfilled.

If $b^r_n=n^{\delta}$, $\delta>0$, then condition \eqref{eq25} is a consequence of \eqref{eq24}.
\end{rem}

To obtain the usual (i.e. a more general than \eqref{HaRe1m}) H\'ajek-R\'enyi type inequality for the moments,
we need the following definition of the Kolmogorov inequality.

\begin{defi}   \label{defKolm2m}
We say that the random variables $X_1, \dots, X_n$ satisfy the second Kolmogorov type maximal inequality for moments,
if for each fixed $k, m$ with $1\leq k \leq  m \leq n$ 
\begin{equation}
\label{Kolm2m}
\ee \Big[ \max_{k \leq l \leq m} |X_k+ \dots + X_l| \Big]^r
\leq
K \sum\nolimits_{l=k}^m \alpha_l
\end{equation}
where $\alpha_1, \dots,\alpha_n$ are non-negative numbers, $r>0$, and $K>0$.
\end{defi}

\begin{defi}   \label{defHaRe2m}
We say that the random variables $X_1, \dots, X_n$ satisfy the second H\'ajek-R\'enyi type maximal inequality for moments,
if for each fixed $m$ with $1\leq m \leq n$  
\begin{equation}
\label{HaRe2m}
\ee \left[ \max_{m\leq l\leq n} \left|\frac{S_l}{\beta_l}\right| \right]^r
\leq
C \left[ \frac{1}{\beta_m^r} \sum_{l=1}^m \alpha_l +  \sum_{l=m+1}^n \frac{\alpha_l}{\beta_l^r} \right] 
\end{equation}
where $\beta_1, \dots,\beta_n$ is a non-decreasing sequence of positive numbers,  
$\alpha_1, \dots,\alpha_n$ are non-negative numbers, 
$r>0$, and $C>0$.
\end{defi}

\begin{thm}
\label{Kolm2m-HaRe2m}
Let the  random variables $X_1, \dots, X_n$ be fixed.
If the second Kolmogorov type maximal inequality \eqref{Kolm2m} for moments is satisfied, 
then the second H\'ajek-R\'enyi type maximal inequality \eqref{HaRe2m} for moments is satisfied with $C= 4 D_r K$,
where $D_r= 1$ for $0<r\leq 1$, and $D_r= 2^{r-1}$ for $r \geq1$.
\end{thm}
\begin{pf}
By the $C_r$-inequality and \eqref{HaRe1m}, 
$$
\ee \left[ \max_{m\leq l\leq n} \left|\frac{S_l}{\beta_l}\right| \right]^r \leq 
D_r \left[ \ee \left( \frac{|S_m|}{\beta_m} \right)^r  + \ee \left( \max_{m+1 \leq l\leq n} \left| \frac{X_{m+1} + \dots +X_l}{\beta_l}\right| \right)^r  \right] \leq 
$$
$$
\leq
D_r \left[ \frac{K}{\beta_m^r} \sum_{l=1}^m \alpha_l +  4 K \sum_{l=m+1}^n \frac{\alpha_l}{\beta_l^r} \right].
$$
\qed\end{pf}
A popular way to obtain the SLLN is the following application of the second H\'ajek-R\'enyi type maximal inequality.

{\it Second proof of Theorem \ref{Kolm1m-SLLN}} if the second Kolmogorov type maximal inequality \eqref{Kolm2m} for moments is satisfied.
 By Theorem~\ref{Kolm2m-HaRe2m},
$$
\ee \left[ \sup_{k\ge m} \left|\frac{S_k}{b_k}\right| \right]^r =
\lim_{n\to\infty} \ee \left[ \max_{m \le  k \le n} \left|\frac{S_k}{b_k}\right| \right]^r \le
C \left[ \frac{1}{b_m^r} \sum_{l=1}^m \alpha_l +  \sum_{l=m+1}^\infty \frac{\alpha_l}{b_l^r} \right]. 
$$
By the Kronecker lemma, the above expression converges to $0$, as $m\to\infty$.
\qed

\section{Applications of the moment inequalities}
\setcounter{equation}{0}

\noindent
{\bf Some basic inequalities}

\noindent
A Kolmogorov type inequality is the starting point of each application of our main theorem.

Doob's inequalities for martingales (in particular for independent random variables) are well-known.

Qi-Man Shao (2000) proved a comparison theorem for moment inequalities between negatively associated and independent random variables. 
The family of random variables $X_1, X_2, \dots$ is called {\it negatively associated}
if for any pair of finite disjoint subsets $A_1$ and $A_2$ of $\{ 1, 2, \dots\}$ 
$$
\cov \left( f_1(X_i, i\in A_1), f_2(X_i, i\in A_2) \right) \le 0
$$
whenever $f_1$ and $f_2$ are coordinate-wise increasing and the above covariance exists.
By Qi-Man Shao (2000), Theorem 1, we have the following.
If the family of random variables $X_1, X_2, \dots$ is negatively associated, 
the random variables $X^*_1, X^*_2, \dots$ are independent such that $X^*_i$ and $X_i$ have the same distribution for each $i=1,\dots, n$,
and $f$ is a non-decreasing convex real function,
then
\begin{equation}
\label{Shao1}
\ee f \Big( \max_{1 \leq k \leq n} |X_1+ \dots + X_k| \Big)
\leq
\ee f \Big( \max_{1 \leq k \leq n} |X^*_1+ \dots + X^*_k| \Big)
\end{equation}
whenever the right hand side exists.

The general inequality \eqref{Shao1} implies the following Kolmogorov type inequalities for negatively associated zero mean random variables $X_1, \dots, X_n$
(see Qi-Man Shao (2000), Theorem 2).
\begin{equation}
\label{Shao2}
\ee \Big[ \max_{1\leq k \leq n} \big|\sum\limits_{i=1}^k X_i \big| \Big]^p
\leq
2^{3-p} \sum\limits_{i=1}^n \ee | X_i |^p  \quad {\text {for}} \quad  1 < p \leq 2,
\end{equation}
\begin{equation}
\label{Shao3}
\ee \Big[ \max_{1\leq k \leq n} \big|\sum\limits_{i=1}^k X_i \big| \Big]^p
\leq
2 \left( \frac{15p}{ \ln p} \right)^{p} \left\{ \left( \sum\limits_{i=1}^n \ee | X_i |^2 \right)^{p/2} +  \sum\limits_{i=1}^n \ee | X_i |^p \right\}  
\quad {\text {for}} \quad  p > 2.
\end{equation}
(See also Su, Zhao, and Wang (1997).)
Inequality \eqref{Shao2} was obtained by Matu{\l}a  (1992) for $p=2$. 

Kuczmaszewska (2005) obtained for negatively associated zero mean random variables $X_1, \dots, X_n$ the following inequality.
\begin{equation}
\label{Kuczma}
\ee \Big[ \max_{1\leq k \leq n} \big|\sum\limits_{i=1}^k X_i \big| \Big]^4
\leq
\sum\limits_{i=1}^n \ee X_i^4  + 2 \sum\limits_{i=1}^n \ee X_i^2 \sum\limits_{j=1}^{i-1} \ee X_j^2 .
\end{equation}

Let $\{X_n,n\ge1\}$ be a sequence of random variables. 
Let $\rho(n)$ be its {\it Kolmogorov--Rozanov mixing coefficient}, that is,
$$
\varrho(n)=\sup_{\scriptstyle{\genfrac{}{}{0pt}{}{k\ge 1}{X\in L_1^k, Y\in L_{k+n}^\infty}}}
\frac {\operatorname{cov} (X,Y)}{\sqrt{\operatorname{var} (X)\cdot
\operatorname{var} (Y)}},
$$
where $L_i^j$ is the space of square integrable random variables that are
measurable with respect to the $\sigma$-field $\mathfrak F_i^j$ generated by $X_i,\dots, X_j$. 
In the case $j=\infty$ the $\sigma$-field $\mathfrak F_i^\infty$ is supposed to be generated by the random variables $X_i,X_{i+1},\dots$ .
A sequence of random variables is said to be $\varrho$-{\it mixing} if its Kolmogorov--Rozanov mixing coefficient $\varrho(n)$ tends to zero, as $n\to\infty$.

Using $\varrho(n)$, Shao (1995) proved the following result:
{if $q\ge2$ and $\E X_n=0$\rom, $\E |X_n|^q<\infty$\rom, $n\ge1$\rom, then there exists a constant $\theta$ depending only on $q$ and $\{\varrho(n)\}$ such that}
\begin{equation}  \label{eq61}
\begin{aligned}
\ee \left[ \max_{k\le n} |S_k| \right]^q\le\theta  & \left[ n^{q/2}\exp \bigg\{\theta\sum_{i=0}^{[\log n]}\varrho(2^i)\bigg\}\max_{j\le n} \left(\ee X_j^2\right)^{q/2}\right.
\\
&\qquad +\left.  n\exp\bigg\{\theta\sum_{i=0}^{[\log n]}\varrho^{2/q}(2^i)\bigg\}\max_{j\le n} \ee |X_j|^q \right].
\end{aligned}
\end{equation}

The simplest case of this inequality is presented by $q=2$.
Indeed, for $q=2$ we have
\begin{equation} \label{eq62}
\ee \left[ \max_{k\le n}|S_k| \right]^2 \le 2\theta n \exp\bigg\{\theta\sum\nolimits_{i=0}^{[\log n]}\varrho(2^i)\bigg\}\max_{j\le n} \ee X_j^2.
\end{equation}
   
\noindent
{\bf Proofs of SLLN's}

\noindent
There are several methods to obtain SLLN's.
Our method described in the previous section is known for specialists.
However, it is worth to fix explicitly the conditions like we did in Theorem~\ref{Kolm1m-SLLN} because it helps to handle some difficult particular cases.

In Fazekas-Klesov (2000) the general SLLN was applied among others to prove Brunk-Prokhorov (see Brunk (1948), Prokhorov (1950)) type SLLN for martingales (Corollary 3.1),
to extend Shao's (1995) Marcinkiewicz-Zygmund type SLLN for $\varrho$-mixing sequences (Theorems 5.1 and 5.2), 
and to extend Hansen's (1991) SLLN for mixingales (Theorems 6.1 and 6.2).

Using the general Theorem~\ref{Kolm1m-SLLN}, Kuczmaszewska (2005) obtained SLLN's for certain dependent random variables.
She presented an SLLN for negatively associated sequences.
She obtained the Kolmogorov type maximal inequality which serves the base of the proof (see \eqref{Kuczma} above).
Moreover, she proved a general Marcinkiewicz-Zygmund type SLLN for $\varrho$-mixing sequences.
Both Kuczmaszewska (2005) and Fazekas-Klesov (2000) applied the Kolmogorov type maximal inequality for  $\varrho$-mixing sequences given in Shao (1995), see
\eqref{eq61} above.

\noindent
{\bf Rate of convergence in the SLLN}

\noindent
Hu and Hu (2006) obtained some results concerning the rate of convergence in the SLLN.
They followed the approach described in Fazekas-Klesov (2000) but they utilized the full strength of Dini's theorem.

\begin{lem}
\label{LemmaDini}
(Dini's theorem, see Fikhtengolts (1969), sect. 375.5.)
Let $c_1, c_2, \dots$ be non-negative numbers, $\nu_n = \sum_{k=n}^\infty c_k$. 
If $0< \nu_n < \infty $ for all $n= 1,2, \dots $, then for any $0 < \delta < 1$ 
we have
$ \sum_{n=1}^\infty  {c_n}/{\nu_n^\delta}  < \infty$.
\end{lem}

\begin{thm}
\label{thmShuheMing} (Lemma  1.2 in Hu-Hu (2006).)
Let $X_1, X_2, \dots$ be a sequence of random variables.
Let the non-negative numbers $\alpha_1, \alpha_2, \dots$, $r>0$, and $K>0$ be fixed.
Assume that for each $m\geq 1$ the first Kolmogorov type maximal inequality \eqref{Kolm1m} for moments is satisfied.
Let $b_1, b_2, \dots$ be a non-decreasing unbounded sequence of positive numbers. 
If \eqref{sum} is satisfied, i.e. 
$
\sum\nolimits_{l=1}^{\infty} \frac{\alpha_l}{b_l^r} < \infty
$,
then
$
\lim_{n\to\infty} \frac{S_n}{b_n} = 0
$ 
a.s., moreover
$$
\frac{S_n}{b_n} = {\rm O} \left( \frac{\beta_n}{b_n} \right) \quad {\text a.s.}
$$
 where 
$$
\beta_n = \max_{1 \leq k \leq n} b_k \nu_k^{\delta/r}, \quad  \nu_k = \sum_{l=k}^\infty \frac{\alpha_l}{b_l^r},
$$
$\delta $ is an arbitrary number with $0 < \delta <  1$. 
\end{thm}
We remark that $\lim_{n\to\infty} \frac{\beta_n}{b_n} = 0$ in the above theorem.

This theorem was used to obtain convergence rates in SLLN's for random variables satisfying certain dependence conditions.
We list those cases where the SLLN was obtained in Fazekas-Klesov (2000) while the rate of convergence in Hu-Hu (2006).
(a) General SLLN for sequences with superadditive moment function (sect. 7 in Fazekas-Klesov (2000) and Theorem 2.1 in Hu-Hu (2006)).
The proofs are based on an inequality in M\'oricz (1976).
(b) Marcinkiewicz-Zygmund SLLN for sequences with superadditive moment function (Theorem 7.1  in Fazekas-Klesov (2000) and Theorem 2.2 in Hu-Hu (2006)).
(c) SLLN using Petrov's natural characteristics (see Petrov (1975)) of the order of growth of sums (Theorem 4.1 in Fazekas-Klesov (2000) and Theorem 2.4 in Hu-Hu (2006)).

Simple versions of the so called almost sure limit theorems are based on an SLLN with logarithmic normalizing factors (see, e.g., M\'ori (1993)).
In  Fazekas-Klesov (2000), Theorem 8.1, a short proof is given for the next SLLN.
Following the lines of that proof, Hu and Hu (2006) gave the rate of convergence.

\begin{thm} 
\label{Theorem91} (Hu-Hu (2006), Theorem 2.5.)
For some $\beta > 0$ and $C>0$ let
\begin{equation} \label{eq91}
| \cov (X_k, X_l) | \leq C \left( \frac lk \right)^{\beta}, \qquad 1 \leq l \leq k .
\end{equation}
Then
\begin{equation} \label{eq92}
\lim_{n\to\infty} \frac{1}{\log{n}} \sum_{k=1}^n \frac{X_k - \ee X_k}{k} = 0 \qquad \text {a.s.}
\end{equation}
Moreover, for any $0 < \delta < 1/2$
\begin{equation} \label{eq92k}
 \frac{1}{\log n} \sum_{k=1}^n \frac{X_k - \ee X_k}{k} = {\rm O} \left( \frac {1}{(\log{n})^\delta } \right)  \qquad \text {a.s.}
\end{equation}
\end{thm}

\section{H\'ajek-R\'enyi type maximal inequalities for the probability}
\setcounter{equation}{0}
%
We start with the abstract Kolmogorov type and H\'ajek-R\'enyi type maximal inequalities for the probability.
\begin{defi}   \label{defKolm1p}
We say that the random variables $X_1, \dots, X_n$ satisfy the first Kolmogorov type maximal inequality for the probability,
if for each $m$ with $1\leq m \leq n$ 
\begin{equation}
\label{Kolm1p}
\pp \Big( \max_{1\leq l\leq m} |S_l| \ge \varepsilon \Big)
\leq
\frac{K}{\varepsilon^r} \sum_{l=1}^m \alpha_l \quad {\text{ for \ any}} \quad \varepsilon >0
\end{equation}
where $\alpha_1, \dots,\alpha_n$ are non-negative numbers, $r>0$, and $K>0$.
\end{defi}

\begin{defi}   \label{defHaRe1p}
We say that the random variables $X_1, \dots, X_n$ satisfy the first H\'ajek-R\'enyi type maximal inequality for the probability,
if 
\begin{equation}
\label{HaRe1p}
\pp \left( \max_{1\leq l\leq n} \left| \frac{S_l}{\beta_l} \right| \ge \varepsilon \right)
\leq
\frac{C}{\varepsilon^r} \sum_{l=1}^n \frac{\alpha_l}{\beta_l^r} \quad {\text{ for \ any}} \quad \varepsilon >0
\end{equation}
where $\beta_1, \dots,\beta_n$ is a non-decreasing sequence of positive numbers,  
$\alpha_1, \dots,\alpha_n$ are non-negative numbers, 
$r>0$, and $C>0$.
\end{defi}

\begin{thm}
\label{Kolm1p-HaRe1p}
(T\'om\'acs-L\'{i}bor (2006), Theorem 2.1.)
Let the  random variables $X_1, \dots, X_n$ be fixed.
If the first Kolmogorov type maximal inequality \eqref{Kolm1p} for the probability is satisfied, 
then the first H\'ajek-R\'enyi type maximal inequality \eqref{HaRe1p} for the probability is satisfied with $C= 4K$.
\end{thm}

The above result allows us to obtain an abstract form of the strong law of large numbers.
\begin{thm}
\label{Kolm1p-SLLN} (Theorem 2.4 in T\'om\'acs-L\'{i}bor (2006).)
Let $X_1, X_2, \dots$ be a sequence of random variables.
Let the non-negative numbers $\alpha_1, \alpha_2, \dots$, $r>0$, and $K>0$ be fixed.
Assume that for each $m\geq 1$ the first Kolmogorov type maximal inequality \eqref{Kolm1p} for the probability is satisfied.
Let $b_1, b_2, \dots$ be a non-decreasing unbounded sequence of positive numbers. 
If
\begin{equation}
\label{sum2}
\sum\nolimits_{l=1}^{\infty} \frac{\alpha_l}{b_l^r} < \infty
\end{equation}
then
\begin{equation}
\label{lim2}
\lim_{n\to\infty} \frac{S_n}{b_n} = 0 \ \ \text {a.s.}
\end{equation}
\end{thm}

To obtain the usual (i.e.  of the shape of \eqref{HaRe}) H\'ajek-R\'enyi type inequality for the probability,
we need the following definition of the Kolmogorov inequality.

\begin{defi}   \label{defKolm2p}
We say that the random variables $X_1, \dots, X_n$ satisfy the second Kolmogorov type maximal inequality for the probability,
if for each fixed $k, m$ with $1\leq k \leq  m \leq n$ 
\begin{equation}
\label{Kolm2p}
\pp \Big( \max_{k \leq l \leq m} |X_k+ \dots + X_l| \ge \varepsilon \Big)
\leq
\frac{K}{\varepsilon^r} \sum_{l=k}^m \alpha_l  \quad {\text{ for \ each}} \quad \varepsilon > 0
\end{equation}
where $\alpha_1, \dots,\alpha_n$ are non-negative numbers, $r>0$, and $K>0$.
\end{defi}

The shape of the following inequality is the same as that of the original H\'ajek-R\'enyi type inequality \eqref{HaRe}.
\begin{defi}   \label{defHaRe2p}
We say that the random variables $X_1, \dots, X_n$ satisfy the second H\'ajek-R\'enyi type maximal inequality for the probability,
if for each fixed $m$ with $1\leq m \leq n$  
\begin{equation}
\label{HaRe2p}
\pp \left( \max_{m\leq l\leq n} \left|\frac{S_l}{\beta_l}\right| \ge \varepsilon \right)
\leq
\frac{C}{\varepsilon^r}  \left[ \frac{1}{\beta_m^r} \sum_{l=1}^m \alpha_l +  \sum_{l=m+1}^n \frac{\alpha_l}{\beta_l^r} \right]  \ \ {\text{ for \ any}} \ \ \varepsilon > 0
\end{equation}
where $\beta_1, \dots,\beta_n$ is a non-decreasing sequence of positive numbers,  
$\alpha_1, \dots,\alpha_n$ are non-negative numbers, 
$r>0$, and $C>0$.
\end{defi}

\begin{thm}
\label{Kolm2p-HaRe2p}
Let the  random variables $X_1, \dots, X_n$ be fixed.
If the second Kolmogorov type maximal inequality \eqref{Kolm2p} for the probability is satisfied, 
then the second H\'ajek-R\'enyi type maximal inequality \eqref{HaRe2p} for the probability is satisfied with $C= (1+ \sqrt[r]{4})^r K$.
\end{thm}
\begin{pf}
Let $0< a < 1$. Then, by \eqref{HaRe1p}, 
$$
\pp \left( \max_{m\leq l\leq n} \left|\frac{S_l}{\beta_l}\right| \ge \varepsilon \right)
\leq
\pp \left( \frac{|S_m|}{\beta_m} \ge a \varepsilon \right)  + 
\pp \left( \max_{m+1 \leq l\leq n} \left| \frac{X_{m+1} + \dots +X_l}{\beta_l}\right| \ge (1-a) \varepsilon \right)  \leq 
$$
$$
\leq
\frac{K}{\varepsilon^r a^r}  \frac{1}{\beta_m^r} \sum_{l=1}^m \alpha_l +  \frac{4 K}{\varepsilon^r (1-a)^r} \sum_{l=m+1}^n \frac{\alpha_l}{\beta_l^r}.
$$
\qed\end{pf}
We show that by applying the second H\'ajek-R\'enyi type maximal inequality we can prove the SLLN.

{\it Second proof of Theorem \ref{Kolm1p-SLLN}} if the second Kolmogorov type maximal inequality for the probability is satisfied.
By Theorem~\ref{Kolm2p-HaRe2p},
$$
\pp \left( \sup_{k\ge m} \left|\frac{S_k}{\beta_k}\right| > \varepsilon \right) 
\le
\frac{C}{\varepsilon^r} \left[ \frac{1}{b_m^r} \sum_{l=1}^m \alpha_l +  \sum_{l=m+1}^\infty \frac{\alpha_l}{b_l^r} \right]. 
$$
As $m\to\infty$, the above expression converges to $0$.
\qed

Finally, we remark that, by Markov's inequality, a Kolmogorov (or a H\'ajek-R\'enyi) type inequality for moments always implies an appropriate inequality for probability.

\section{Applications of the probability inequalities}
\setcounter{equation}{0}

\noindent
{\bf Some classical H\'ajek-R\'enyi type maximal inequalities}

\noindent
Most of the classical H\'ajek-R\'enyi type inequalities were proved by direct methods.
Below we list some known results and point out how can we insert them into our general framework.
However, by our general method, we can reproduce the H\'ajek-R\'enyi type inequalities up to an absolute constant multiplier. 

First we remark that in H\'ajek-R\'enyi (1955) a direct proof is given for inequality \eqref{HaRe} if the random variables are independent.
Now we see that the original H\'ajek-R\'enyi inequality is a consequence of the original Kolmogorov inequality (up to constant multiplier $C=9$).

Another classical H\'ajek-R\'enyi type inequality was proved by Chow (1960)  for submartingales.
We can see that it is a consequence of Doob's inequality (up to a constant).

Kounias and Weng (1969) proved the following  H\'ajek-R\'enyi type inequality for arbitrary random variables (i.e. without assuming any dependence condition).
Let $X_1, \dots, X_n$ be random variables. 
Let $r>0$ be fixed. 
Let the moments $v_i = \ee |X_i|^r$ be finite for each $i$.
Let $s=1$ if $0<r\le 1$ and $s=r$ if $r>1$.
Then
\begin{equation}
\label{KouWeng}
\pp \left( \max_{1\leq l\leq n} \left| \frac{S_l}{\beta_l} \right| \ge \varepsilon \right)
\leq
\frac{1}{\varepsilon^r} \left( \sum_{l=1}^n \Big( \frac{v_l}{\beta_l^r} \Big)^{1/s} \right)^s   \quad {\text{ for \ any}} \quad \varepsilon >0
\end{equation}
where $\beta_1, \dots,\beta_n$ is an arbitrary  non-decreasing sequence of positive numbers. 

This inequality is of type \eqref{HaRe1p} for $r \leq 1$.
But for $r>1$ it seems to be different of type \eqref{HaRe1p}.
However, we shall see that it can be inserted into our framework for $r>1$, too.
The appropriate Kolmogorov type inequality is of the form \eqref{Kolm1p}
with 
$
\alpha_k = \left( \sum_{i=1}^k a_i \right)^r - \left( \sum_{i=1}^{k-1} a_i \right)^r
$,
$a_i= (\ee |X_i|^r)^{1/r}$.
We can prove that 
$
\sum_{l=1}^n \frac{\alpha_l}{\beta_l^r} \leq  \left( \sum_{l=1}^n \frac{a_l}{\beta_l} \right)^r
$.
So in this case \eqref{HaRe1p} implies \eqref{KouWeng} up to a multiplier 4.

Szynal (1973) obtained the following  H\'ajek-R\'enyi type inequality for arbitrary random variables (i.e. without assuming independence and moment conditions).
Let $X_1, \dots, X_n$ be arbitrary random variables. 
Let $r>0$ be fixed. 
Let $s=1$ if $0<r\le 1$ and $s=r$ if $r>1$.
Then
\begin{equation}
\label{Szynal}
\pp \left[ \max_{m \leq l\leq n} \left| \frac{S_l}{\beta_l} \right| \ge 3 \varepsilon \right]
\leq
 2 \left[ \sum_{l=1}^m    \ee^{1/s }\Big( \frac{|X_l|^r}{ ( \beta_m \varepsilon )^r +|X_l|^r} \Big) 
+         \sum_{l=m+1}^n  \ee^{1/s }\Big( \frac{|X_l|^r}{ ( \beta_l \varepsilon )^r +|X_l|^r} \Big) \right]^s  
\end{equation}
for any $\varepsilon >0$ where $\beta_1, \dots,\beta_n$ is an arbitrary  non-decreasing sequence of positive numbers. 
It seems that \eqref{Szynal} can not be inserted into the framework of our main theorem.

Bickel (1970) obtained a H\'ajek-R\'enyi type generalization of L\'evy's inequality.

\noindent
{\bf Rate of convergence in SLLN's}

\noindent
Combining the method of Theorem~\ref{Kolm1p-SLLN} and the ideas of Hu and Hu (2006), we can obtain further results concerning the rate of convergence in the SLLN.

\begin{thm}
\label{Libor} 
Let $X_1, X_2, \dots$ be a sequence of random variables.
Let the non-negative numbers $\alpha_1, \alpha_2, \dots$, $r>0$, and $K>0$ be fixed.
Assume that for each $m\geq 1$ the first Kolmogorov type maximal inequality \eqref{Kolm1p} for the probability is satisfied.
Let $b_1, b_2, \dots$ be a non-decreasing unbounded sequence of positive numbers. 
If \eqref{sum2}   is satisfied, i.e. 
$
\sum\nolimits_{l=1}^{\infty} \frac{\alpha_l}{b_l^r} < \infty
$,
then
$
\lim_{n\to\infty} \frac{S_n}{b_n} = 0
$ 
a.s., moreover
$$
\frac{S_n}{b_n} = {\rm O} \left( \frac{\beta_n}{b_n} \right) \quad {\text a.s.}
$$
 where $\beta_n$ is defined in Theorem~\ref{thmShuheMing}. 
\end{thm}
This theorem can be used to obtain convergence rates in the SLLN for random variables satisfying certain dependence conditions.

\noindent
{\bf Alternative proofs for SLLN's and rates of convergences}

Here we show that the general theorems on the SLLN's and on the rates of convergences in them can be applied for sequences with certain dependence properties.

Chandra and Ghosal (1996) obtained a Kolmogorov type inequality and a Marcinkiewicz-Zygmund type SLLN for asymptotically almost negatively associated (AANA) random variables.

The sequence $X_1, X_2, \dots $ is called AANA if there exists a non-negative sequence $q_n \to 0$
such that
\begin{equation}
\label{AANA}
\cov \left( f(X_m), g(X_{m+1}, \dots , X_{m+k}) \right) \leq  q_m  \left(  \var\big(  f(X_m) \big)  \var\big(  g(X_{m+1}, \dots , X_{m+k}) \big) \right)^{1/2}
\end{equation}
for all $m,k \ge 1$ and for all coordinatewise increasing continuous functions $f$ and $g$ whenever the right hand side of the above inequality is finite.

The negatively associated (in particular the independent) sequences are AANA.
The Kolmogorov type inequality for AANA sequences is the following.

\begin{thm} (Chandra-Ghosal (1996), Theorem 1.) Let $X_1, \dots, X_n$ be zero mean square integrable r.v.'s such that \eqref{AANA} holds
for every $1 \leq m < m+k \leq n$. 
Let $A_n = \sum_{l=1}^{n-1} q_l^2$.
Then 
\begin{equation}
\label{ChaGho}
\pp \Big( \max_{1\leq l\leq n} |S_l| \ge \varepsilon \Big)
\leq
\frac{2}{\varepsilon^2}  \left( A_n + \sqrt{1+A_n^2 } \right)^2 \sum_{l=1}^n  \ee X^2_l \quad {\text{ for \ any}} \quad \varepsilon >0 .
\end{equation}
\end{thm}
Using \eqref{ChaGho}, Kim, Ko and Lee (2004) obtained H\'ajek-R\'enyi inequalities of the form 
\eqref{HaRe1p} and \eqref{HaRe2p}, moreover, if $A = \sum_{l=1}^{\infty} q_l^2 < \infty$, an SLLN of the form of Theorem~\ref{Kolm1p-SLLN}.
Now we can see that, for AANA sequences, our method gives not only the SLLN but also the rate of convergence in it.

We mention that negatively associated sequences are AANA. 
In the special case of negatively associated sequences, the Kolmogorov inequality was obtained by Matu{\l}a (1992), 
the H\'ajek-R\'enyi type inequality and the SLLN by Liu, Gan and Chen (1999), while the rate of convergence in the SLLN by Hu and Hu (2006). 

Christofides (2000) obtained a  H\'ajek-R\'enyi type inequality for demisubmartingales.
Let $S_1, S_2, \dots $ be a sequence of random variables with finite expectation.
If
\begin{equation}
\label{demimart}
\ee [(S_{j+1} - S_j ) f(S_1, \dots, S_j ) ] \ge 0
\end{equation}
for any coordinatewise non-decreasing function $f$ such that the above expectation is defined, then $\{ S_j \}_{j \ge 1} $ is called a {\it demimartingale}.
If, in addition, the function $f$ in \eqref{demimart} is assumed to be non-negative, then the sequence $\{ S_j \}_{j \ge 1} $ is called a {\it demisubmartingale}
(see Newman and Wright (1982)). 
We remark that if we omit that the function $f$ is non-decreasing, then \eqref{demimart} is equivalent that $\{ S_j \}_{j \ge 1} $ is a martingale. 

Christofides (2000), Theorem 2.1 is the following H\'ajek-R\'enyi type inequality for demisubmartingales.
Let $S_0, S_1, S_2, \dots $ be a demisubmartingale with $S_0 = 0 $, $\beta_1, \beta_2, \dots $ a non-decreasing sequence of positive numbers, then 
\begin{equation}
\label{Christo}
\pp \left( \max_{1\leq l\leq n} \left| \frac{S_l}{\beta_l} \right| \ge \varepsilon \right)
\leq
\frac{1}{\varepsilon} \sum_{l=1}^n \frac{\ee( S_l^+ - S_{l-1}^+ )} {\beta_l} \quad {\text{ for \ any}} \quad \varepsilon >0
\end{equation}
where $S_l^+$ denotes the positive part of $S_l$.

Using this inequality, an SLLN can be obtained for demimartingales (Christofides (2000), Theorem 2.2).
Now we can see that our general theorems immediately imply the SLLN from the Kolmogorov type inequality, moreover the rate of convergence in the SLLN.

The sequence of partial sums of zero mean associated random variables is a demimartingale.
We mention that for associated random variables  Kolmogorov type inequalities were obtained by Newman and Wright (1982) and by Matu{\l}a (1996).
From those inequalities  Matu{\l}a (1996) and Hu-Hu (2006) obtained the SLLN and the rate of convergence in the SLLN, respectively.
We see that the later theorems are covered by our general method.

\section{A Marcinkiewicz-Zygmund type SLLN}
\setcounter{equation}{0}

In this section we present an abstract Marcinkiewicz-Zygmund type SLLN, i.e. an SLLN without any dependence condition.
We shall use the following notation.
Let $X_1, X_2, \dots$ be a sequence of random variables.
Let
$$
Y_k = X_k {\ii} \{ |X_k| \le k^{1/p} \}+ k^{1/p} {\ii} \{ X_k > k^{1/p} \} -  k^{1/p} {\ii} \{ X_k < - k^{1/p} \},
$$
$$
Z_k = X_k {\ii} \{ |X_k| \le k^{1/p} \},
$$
denote the truncated variables.
Moreover, let
$$
G(y)= \sup_{n\ge 1} \frac{1}{n} \sum\nolimits_{k=1}^n \pp ( |X_k| > y ), \quad  y \ge 0.
$$
Consider the conditions
\begin{equation} \label{pMoment}
\int_{0}^\infty y^{p-1} G(y) dy < \infty, \quad \sum_{k=1}^\infty  \pp ( |X_k|^p > k ) < \infty.
\end{equation}  

\begin{thm}
\label{MarZyg-SLLN} 
Let $X_1, X_2, \dots$ be a sequence of random variables, $S_n = X_1 + \dots + X_n$, $n=1, 2, \dots$, $S_0 = 0$.
Let $0< p < 2$.
Assume that condition \eqref{pMoment} is satisfied.


Moreover, for $1 \le  p < 2 $, assume that  $\ee X_k =0 $ $(k=1,2,\dots)$ and there exists a $K < \infty$ such that for any $n \ge 1$
\begin{equation}
\label{KolmMZ}
\pp \Big( \max_{1 \leq k \leq n} |(Y_1 - \ee Y_1) + \dots + (Y_k - \ee Y_k) | \ge \varepsilon \Big)
\leq
\frac{K}{\varepsilon^2} \sum_{i=1}^n  \dd^2 Y_i  \quad {\text{ for \ any}} \quad \varepsilon >0.
\end{equation}
Then
\begin{equation}
\label{limMZ}
\lim_{n\to\infty} \frac{S_n}{n^{1/p}} = 0 \ \ \text {a.s.}
\end{equation}
Relation \eqref{limMZ} remains valid if in \eqref{KolmMZ} the sequence $\{ Y_k \}$ is replaced by $\{ Z_k \}$.  
\end{thm}

\begin{pf}
The proof of this theorem is similar to that of Theorem 2 in Chandra-Ghosal (1996).
\qed\end{pf}

\section{Further results}
\label{uj}
\setcounter{equation}{0}

\noindent
{\bf General strong laws of large numbers}

\noindent
First we mention a well-known result of Serfling (see Theorem 2.4.1 in Stout (1974)).
Let $X_1, \dots , X_n$ be random variables and let $F_{a,n}$ denote the joint distribution function of $X_{a+1}, \dots , X_{a+n}$.
Assume that $g$ is a real function such that $g(F_{a,k}) + g(F_{a+k,m}) \le g(F_{a,k+m})$ for all $1\le k< k+m$ and $a\ge 0$.
If 
$$
\ee \left( \sum\nolimits_{i=a+1}^{a+n} X_i \right)^2 \le g(F_{a,n}) 
$$
for all $n\ge 1$ and $a\ge 0$, then
\begin{equation} \label{Serfling}
\ee \left( \max_{a<k\le n}\sum\nolimits_{i=a+1}^{a+k} X_i \right)^2 \le \left(\frac{\log(2n)}{\log 2}\right)^2 g(F_{a,n}) 
\end{equation}
for all $n\ge 1$ and $a\ge 0$.

Sung (2008) studies a sequence of random variables without assuming any known weak dependence condition.
The conditions used in the paper are in terms of the covariances of the random variables.
Applying Serfling's inequality \eqref{Serfling}, first a Kolmogorov type maximal inequality for the moments is proved.
Then two general SLLN's are obtained, the second one rests on our Theorem \ref{Kolm1m-SLLN}.

In Shen-Yang-Hu (2013) a general Kolmogorov type inequality for the moments is considered, that is 
\begin{equation} \label{moment2}
\ee \left(\max_{m\le k\le n} \left(\sum\nolimits_{i=m}^{k} X_i \right)^2 \right) \le  K \ee \sum\nolimits_{i=m}^{n} X_i^2 
\end{equation}
for all  $n\ge m \ge 1$, where the constant $K$ is independent of $n$ and $m$.
Two general SLLN's are proved in the case when the truncated and normalized random variables satisfy condition \eqref{moment2}.
Then a Marcinkiewicz SLLN is obtained when the random variables themselves obey \eqref{moment2}.
The proof of the Marcinkiewicz SLLN is based on Theorem \ref{thmShuheMing}.

Yang-Su-Yu (2008) obtained the following general SLLN. 
Let $b_1, b_2, \dots$ be a non-decreasing sequence of positive numbers with $1\le b_{2n} /b_n \le c <\infty$ for some $c>1$.
Assume that for any $\varepsilon>0$,
$$
\sum_{n=1}^\infty \frac{1}{n} \pp \left( \max_{1\le k\le n}  |S_k| > b_n \varepsilon \right) <\infty.
$$
Then 
$$
\lim_{n\to\infty}  \frac{\max_{1\le k\le n} |S_k| }{b_n}  = 0 \quad {\text{a.s.}}
$$
The authors claim that their conditions in certain interesting cases are more general than the ones in Theorem \ref{Kolm1m-SLLN}.

\noindent
{\bf Results for weakly dependent sequences}

Kuczmaszewska (2008) obtained Chung-Teicher type SLLN's for $\varrho^*$-mixing (in other words $\tilde\varrho$-mixing) random variables.
She applied the Rosenthal type inequality given by Utev-Peligrad (2003) and the approach suggested by Fazekas-Klesov (2000).

Wang-Hu-Shen-Ling (2008)  studied strongly positive dependent, (positively) associated, $\tilde\varphi$-mixing, $\tilde\varrho$-mixing, 
and  pairwise negatively quadrant dependent sequences. 
They started with known (Kolmogorov or Rosenthal type) maximal inequalities for the moments and  applying the method of Fazekas-Klesov (2000), Fazekas (2006), 
and Hu-Hu (2006), obtained strong laws of large numbers and rates of convergence of shape of Theorem \ref{thmShuheMing}.

Wang-Hu-Shen-Yang (2010) studied $\psi$-mixing sequences. 
Applying an inequality of Shao (1993), they obtained the following Rosenthal type inequality.
Let $X_1,X_2, \dots$ be a $\psi$-mixing sequence of random variables with mixing coefficient $\psi(n)$ satisfying $\sum_{n=1}^\infty \psi(n) < \infty$.
Let $q\ge 2$ and assume that $\ee| X_n|^q <\infty$,  $\ee X_n =0$ for any $n\ge 1$.
Then there exists a constant $c$ depending only on $q$ and $\psi(.)$ such that
$$
\ee \left( \max_{1\le j\le n} \left| \sum\nolimits_{i=1}^j X_i \right|^q \right) \le c 
\left[ \sum\nolimits_{i=1}^n \left|X_i \right|^q + \left( \sum\nolimits_{i=1}^n X_i^2 \right)^{q/2} \right].
$$
Applying the method of Fazekas (2006) and the above inequality with $q=2$, Wang-Hu-Shen-Yang (2010) obtained H\'ajek-R\'enyi type inequalities and SLLN's.

Yang-Shen-Hu-Wang (2012) studied mixingales and pairwise negatively quadrant dependent (NQD) sequences.
First they obtained Kolmogorov type maximal inequalities for moments.
To this end, for NQD sequences, they applied Serfling's inequality  \eqref{Serfling}.
Then they derived H\'ajek-R\'enyi type inequalities for the probabilities and SLLN's similar to Theorem \ref{thmShuheMing}.
To prove the H\'ajek-R\'enyi type inequalities they applied similar calculations as the ones in Fazekas-Klesov (2000).
However, the consideration of the present paper show, that the H\'ajek-R\'enyi type inequalities are always consequences of the appropriate Kolmogorov inequalities.

Hu-Li-Yang-Wang (2011) studied linear, $\varphi$-mixing,
and  linearly negative quadrant dependent (LNQD) sequences. 
First, using the Longnecker-Serfling  (see Longnecker-Serfling (1977))  method, they obtained Kolmogorov type maximal inequalities from known moment inequalities.
From the maximal inequalities then  they derived H\'ajek-R\'enyi type inequalities and SLLN's.
Now, we see that using our general framework, the proofs of the H\'ajek-R\'enyi inequalities are standard.

Wang-Hu-Li-Yang (2011) improved the results of Kim, Ko and Lee (2004) for AANA sequences.
They started with a Kolmogorov type maximal inequality for moments and derived H\'ajek-R\'enyi type inequalities for the probabilities and SLLN's.
We can see that our general method offers  similar H\'ajek-R\'enyi type inequalities.

The random variables $X_1,X_2, \dots, X_n$ are said to be  negatively orthant dependent (NOD), if
$$
\pp(X_i > x_i, i=1,2,\dots, n) \le \prod_{i=1}^n \pp(X_i > x_i)
$$
and
$$ \pp(X_i \le x_i, i=1,2,\dots, n) \le \prod_{i=1}^n \pp(X_i \le x_i).
$$
Kim (2006) derived H\'ajek-R\'enyi inequalities and an SLLN for NOD random variables.
Songlin (2013) improved the constant multiplier in Kim's H\'ajek-R\'enyi type inequality.
Actually, Songlin applied the same method as the one used for the general proof in T\'om\'acs-L\'{i}bor (2006).
Therefore the result is a direct consequence of Theorem \ref{Kolm1p-HaRe1p}.

In Hu-Chen-Wang (2011) Brunk-Prokhorov type strong laws of large numbers are obtained for martingales and for demimartingales. 
The results are extensions of certain results of Fazekas-Klesov (2000).
To prove the theorems the authors apply the abstract H\'ajek-R\'enyi type approach described in the present paper.

Let $S_1, S_2, \dots $ be a sequence of random variables.
If
\begin{equation}
\label{Ndemimart}
\ee [(S_{j+1} - S_j ) f(S_1, \dots, S_j ) ] \le 0
\end{equation}
for any $j=1,2,\dots$, and any coordinate-wise non-decreasing function $f$ such that the above expectation is defined, then $\{ S_j \}_{j \ge 1} $ is called an {\it N-demimartingale}.
In Wang-Hu-Prakasa Rao-Yang (2011) a H\'ajek-R\'enyi type inequality and an SLLN are derided for N-demimartingales  (see also Theorem 3.7.2 in the book Prakasa Rao (2012)).
The authors apply the method offered by Fazekas-Klesov (2000).

\noindent
{\bf Rate of convergence in the SLLN}

\noindent
Sung-Hu-Volodin (2008) improved the results of Hu-Hu (2006).
They started with the following generalization of Dini's theorem.

\begin{lem} \label{LemmaDini2}
(Lemma 1 of Sung-Hu-Volodin (2008).)
Let $\varphi$ be a positive valued function defined on the positive half line satisfying 
\begin{equation} \label{volodin-fi}
\sum_{n=1}^\infty \frac{\varphi(n)}{n^2} < \infty \ \ {\rm and} \ \ 0< \varphi(x) \uparrow \infty \ \ {\rm on} \ \ [c,\infty) \ \   {\rm for \ some} \ \ c>0.
\end{equation}
Let $c_1, c_2, \dots$ be non-negative numbers such that $c_n> 0$ for infinitely many $n$, $\nu_n = \sum_{k=n}^\infty c_k$. 
If $\sum_{k=1}^\infty c_k < \infty $, then 
we have
$ \sum_{n=1}^\infty  {c_n} \varphi(1/{\nu_n})  < \infty$.
\end{lem}

Using Lemma \ref{LemmaDini2} instead of Dini's theorem, Sung-Hu-Volodin (2008) improved the result of Hu-Hu (2006) (that is Theorem \ref{thmShuheMing}).
\begin{thm}
\label{Volodin} (Theorem  1 in Sung-Hu-Volodin (2008).)
Let $X_1, X_2, \dots$ be a sequence of random variables.
Let the non-negative numbers $\alpha_1, \alpha_2, \dots$, $r>0$, and $K>0$ be fixed.
Assume that infinitely many $\alpha_i$ are positive.
Assume that for each $m\geq 1$ the first Kolmogorov type maximal inequality \eqref{Kolm1m} for moments is satisfied for any $m$.
Let $b_1, b_2, \dots$ be a non-decreasing unbounded sequence of positive numbers. 
Assume that \eqref{sum} is satisfied, i.e. 
$
\sum\nolimits_{l=1}^{\infty} \frac{\alpha_l}{b_l^r} < \infty.
$
Let $\varphi$ be a function satisfying \eqref{volodin-fi} and let 
$$
\beta_n = \max_{1 \leq k \leq n} b_k (\varphi(1/\nu_k))^{-1/r}, \quad  {\rm where} \quad \nu_k = \sum_{l=k}^\infty \frac{\alpha_l}{b_l^r}.
$$
Then the following statements hold.
If the sequence $\{\beta_n, \, n\ge 1 \}$ is bounded,   then $S_n/\beta_n$ is almost surely bounded.
If the sequence $\{\beta_n, \, n\ge 1 \}$ is unbounded, then $S_n/\beta_n$ converges almost surely to $0$.
\end{thm}
We remark that $\lim_{n\to\infty} \frac{\beta_n}{b_n} = 0$.
Therefore in the above theorem at each case
$
\lim_{n\to\infty} \frac{S_n}{b_n} = 0
$ 
a.s.
The above theorem with $\varphi(x) = |x|^\delta$ ($0<\delta<1$) gives Theorem \ref{thmShuheMing}.

A version of Theorem \ref{Volodin} was proved by T\'om\'acs (2007) when instead of the moment inequality \eqref{Kolm1m} the probability inequality \eqref{Kolm1p} was assumed.

\end{document}